\documentclass{amsart}

\usepackage{amsfonts}
\usepackage{graphicx}
\usepackage{epic,eepic}

\newcommand{\sech}{\mathop{\rm sech}\nolimits}
\newcommand{\beq}{\begin{equation}}
\newcommand{\eeq}{\end{equation}}

\begin{document}

\title[Lax matrices for Yang-Baxter maps]{Lax matrices for Yang-Baxter maps}
       \author{Yuri Suris}
       \address{$^*$Institut f\"ur Mathematik, Technische Universit\"at
Berlin, Str. des 17. Juni 136, 10623 Berlin, Germany}
           \email{suris@sfb288.math.tu-berlin.de}

       \author{Alexander Veselov}
       \address{Department of Mathematical Sciences,
        Loughborough University, Loughborough,
        Leicestershire, LE11 3TU, UK and
       Landau Institute for Theoretical Physics, Moscow, Russia}
       \email{A.P.Veselov@lboro.ac.uk}

\maketitle




{\small  {\bf Abstract.} It is shown that for a certain class of
Yang-Baxter maps (or set-theoretical solutions to the quantum 
Yang-Baxter equation) the Lax representation can be derived
straight from the map itself. A similar phenomenon for 3D
consistent equations on quad-graphs has been recently discovered
by A. Bobenko and one of the authors, and by F. Nijhoff.}

\bigskip

\subsection*{Introduction}

In 1990 V.G. Drinfeld suggested the problem of studying the
solutions of the quantum Yang-Baxter equation in the case when the
vector space $V$ is replaced by an arbitrary set $X$ and tensor
product by the direct product of the sets (``set-theoretical
solutions to the quantum Yang-Baxter equation'') \cite{D}. In the paper
\cite{V} one of the authors investigated the dynamical aspects of
this problem and suggested a shorter term ``Yang-Baxter map'' for
such solutions.

For each Yang-Baxter map one can introduce the hierarchy of
commuting transfer-maps which are believed to be integrable (see
\cite{V}). In this note we explain how to find Lax representations
for a certain class of Yang-Baxter maps thus giving another
justification for this conjecture. We were motivated by the
explicit examples of the Yang-Baxter maps from \cite{V} and recent
results on the equations on quad-graphs, satisfying the so-called
``3D consistency condition'' \cite{BS, N}.

\subsection*{Yang-Baxter maps and their Lax representations}

Let $X$ be any set and $R$ be a map: $$R: X \times X \rightarrow X
\times X.$$ Let $R_{ij}: X^{n} \rightarrow X^{n}, \quad X^{n} = X
\times X \times .....\times X$ be the map which acts as $R$ on
$i$-th and $j$-th factors and identically on the others. Let
$R_{21}=PRP$, where $P: X^2 \rightarrow X^2$ is the permutation:
$P(x,y) = (y,x)$.

Following \cite{V}, we call $R$ the {\it Yang-Baxter map} if it
satisfies the Yang-Baxter relation
\begin{equation}\label{YB}
R_{23} R_{13} R_{12} = R_{12} R_{13} R_{23},
\end{equation}
considered as the equality of the maps of $X \times X \times X$
into itself.  If additionally $R$ satisfies the relation
\begin{equation}
\label{U} R_{21} R = Id,
\end{equation}
it is called {\it reversible Yang-Baxter map}. Reversibility
condition will not play an essential role in this note but it is
satisfied in all the examples we present.

\begin{figure}[htbp]
\setlength{\unitlength}{0.06em}
\begin{center}
\begin{picture}(550,150)(25,0)
  \put(30,10){\line(5,3){200}}
  \put(30,140){\line(5,-3){200}}
  \put(30,100){\line(1,0){200}}
  \put(10,10){$z$}
  \put(10,95){$y$}
  \put(10,140){$x$}
  \put(125,110){$y_1$}
  \put(100,80){$x_2$}
  \put(165,80){$z_1$}
  \put(240,130){$z_{12}$}
  \put(240,95){$y_{13}$}
  \put(240,15){$x_{23}$}
  \put(350,10){\line(5,3){200}}
  \put(350,140){\line(5,-3){200}}
  \put(350,50){\line(1,0){200}}
  \put(330,10){$z$}
  \put(330,45){$y$}
  \put(330,140){$x$}
  \put(455,35){$y_3$}
  \put(420,70){$z_2$}
  \put(490,70){$x_3$}
  \put(560,130){$z_{12}$}
  \put(560,45){$y_{13}$}
  \put(560,15){$x_{23}$}
  \put(290,70){$=$}
\end{picture}
\end{center}
\caption{Standard representation of the Yang--Baxter relation}
\label{Fig:YB0}
\end{figure}
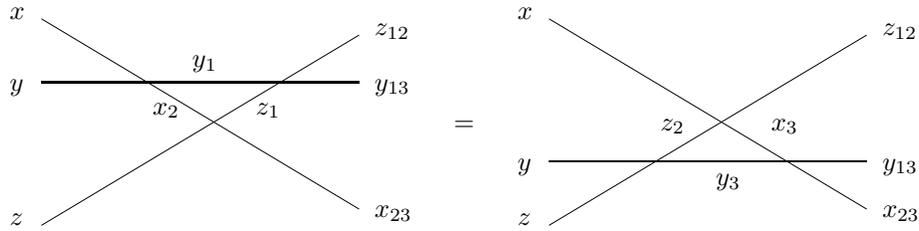

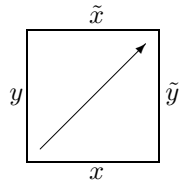
\begin{figure}[bp]
\begin{center}
\setlength{\unitlength}{0.05em}
\begin{picture}(200,140)(-50,-20)
  \put( 0,  0){\line(1,0){100}}
  \put( 0,100){\line(1,0){100}}
  \put(  0, 0){\line(0,1){100}}
  \put(100, 0){\line(0,1){100}}
  \put(10,10){\vector(1,1){80}}
  \put(47,-13){$x$}
  \put(-13,47){$y$}
  \put(47,105){$\tilde{x}$}
  \put(105,47){$\tilde{y}$}
\end{picture}
\caption{A map associated to an elementary quadrilateral; fields
are assigned to edges} \label{Fig: map}
\end{center}
\end{figure}

\begin{figure}[bp]
\setlength{\unitlength}{0.06em}
\begin{center}
\begin{picture}(450,170)(-30,-20)
  \put( 0,  0){\line(1,0){100}}
  \put(100,  0){\line(5,3){50}}
  \put(150,30){\line(0,1){100}}
  \put(50,130){\line(1,0){100}}
  \multiput(50,30)(20,0){5}{\line(1,0){15}}
  \put(  0, 0){\line(0,1){100}}
  \multiput(50,30)(0,20){5}{\line(0,1){15}}
  \put(  0,100){\line(5,3){50}}
  \multiput(50,30)(-16.67,-10){3}{\line(-5,-3){12}}
  \put(40,35){\vector(-1,2){30}}
  \put(140,40){\vector(-1,1){80}}
  \put(95,5){\vector(-2,1){40}}
  \put(105,85){$R_{13}$}
  \put(23,80){$R_{23}$}
  \put(40,5){$R_{12}$}
  \put(40,-11){$x$}
  \put(130,6){$y$}
  \put(155,75){$z$}
  \put(90,135){$x_{23}$}
  \put(5,120){$y_{13}$}
  \put(-25,50){$z_{12}$}
  \put(90,40){$x_2$}
  \put(20,25){$y_1$}
  \put(55,75){$z_1$}
  \put(190,65){$=$}
  \put(250,  0){\line(1,0){100}}
  \put(250, 0){\line(0,1){100}}
  \put(250,100){\line(1,0){100}}
  \put(350, 0){\line(0,1){100}}
  \put(350, 0){\line(5,3){50}}
  \put(400,30){\line(0,1){100}}
  \put(350,100){\line(5,3){50}}
  \put(300,130){\line(1,0){100}}
  \put(250,100){\line(5,3){50}}
  \put(350,100){\line(5,3){50}}
  \put(390,35){\vector(-1,2){30}}
  \put(340,10){\vector(-1,1){80}}
  \put(345,105){\vector(-2,1){40}}
  \put(280,30){$R_{13}$}
  \put(355,45){$R_{23}$}
  \put(290,105){$R_{12}$}
  \put(290,-11){$x$}
  \put(380,6){$y$}
  \put(405,75){$z$}
  \put(340,135){$x_{23}$}
  \put(255,120){$y_{13}$}
  \put(225,50){$z_{12}$}
  \put(300,87){$x_3$}
  \put(368,100){$y_3$}
  \put(332,50){$z_2$}
\end{picture}
\end{center}
\caption{``Cubic'' representation of the Yang--Baxter relation}
\label{Fig:YB}
\end{figure}
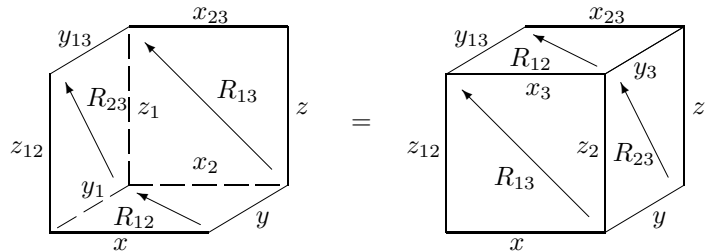

The standard way to represent the Yang-Baxter relation is given by
the diagram in Fig. \ref{Fig:YB0}.
However we would like to use here also an alternative (dual) way to
visualize it, which emphasizes the relation with 3D consistency condition
for discrete equations on quad-graphs (see \cite{BS,ABS}). In
this representation the fields (elements of $X$) are assigned to
the edges of elementary quadrilaterals, so that Fig.\ref{Fig: map}
encodes the map $R:(x,y)\mapsto(\tilde{x},\tilde{y})$.
Then the Yang-Baxter relation is illustrated as  in Fig.
\ref{Fig:YB}.
It encodes the 3D consistency of the maps $R$ attached to all
facets of an elementary cube. The left--hand side of (\ref{YB})
corresponds to the chain of maps along the three rear faces of the
cube on Fig. \ref{Fig:YB}:
\[
R_{12}:(x,y)\mapsto(x_2,y_1),\quad
R_{13}:(x_2,z)\mapsto(x_{23},z_1),\quad
R_{23}:(y_1,z_1)\mapsto(y_{13},z_{12}),
\]
while its right--hand side corresponds to the chain of the maps
along the three front faces of the cube:
\[
R_{23}:(y,z)\mapsto(y_3,z_2),\quad
R_{13}:(x,z_2)\mapsto(x_3,z_{12}),\quad R_{12}:(x_3,y_3)\mapsto
(x_{23},y_{13}).
\]
So, (\ref{YB}) assures that two ways of obtaining
$(x_{23},y_{13},z_{12})$ from the initial data $(x,y,z)$ lead to
the same results.

One can consider also {\it parameter-dependent Yang-Baxter maps}
$R(\lambda,\mu)$, with $\lambda, \mu \in {\bf C}$, satisfying the
corresponding version of Yang-Baxter relation
\begin{equation}\label{sYB}
R_{23}(\mu,\nu) R_{13}(\lambda,\nu) R_{12}(\lambda,\mu) =
R_{12}(\lambda,\mu) R_{13}(\lambda,\nu) R_{23} (\mu,\nu).
\end{equation}
The reversibility condition in this situation reads
\begin{equation}\label{sU}
R_{21}(\mu,\lambda) R(\lambda,\mu) = Id.
\end{equation}
One thinks of the parameters $\lambda,\mu$ as assigned to the same
edges of the quadrilateral in Fig. \ref{Fig: map} as the fields
$x,y$ are. Moreover, opposite edges are thought of as carrying the
same parameters. Thus, in Fig. \ref{Fig:YB} all edges parallel to
the $x$ (resp. $y,z$) axis, carry the parameter $\lambda$ (resp.
$\mu,\nu$). Although this can be considered as a particular case
of the general notion, by introducing $\tilde X = X \times{\bf C}$
and $\tilde R (x,\lambda; y, \mu) = R (\lambda,\mu) (x,y)$, it is
convenient for us to keep the parameter separately.

By the {\it Lax matrix} (or Lax representation) for such a map we
will mean the matrix $A(x,\lambda; \zeta)$ depending on the point
$x \in X$, parameter $\lambda$ and additional (``spectral'')
parameter $\zeta \in {\bf C},$  which satisfies the following
relation:
\begin{equation}\label{Lax}
 A(x,\lambda; \zeta)A(y,\mu; \zeta) =
 A(\tilde{y},\mu; \zeta)A(\tilde{x},\lambda; \zeta),
\end{equation}
whenever $(\tilde{x},\tilde{y}) = R (\lambda,\mu) (x,y).$ As it
was shown in \cite{V}, such a matrix allows one to produce
integrals for the dynamics of the related transfer-maps.

Our main result is the following observation.

Suppose that on the set $X$ we have an action of the linear group
$G =  GL_N$, and that the Yang-Baxter map $R(\lambda, \mu)$ has
the following special form:
\begin{equation}\label{map}
\tilde{x} = B(y,\mu, \lambda)[x], \quad  \tilde{y} =
A(x,\lambda, \mu)[y],
\end{equation}
where $A, B: X \times {\bf C} \times {\bf C} \rightarrow GL_N$ are
some matrix valued functions on $X$ depending on parameters
$\lambda$ and $\mu$ and $A[x]$ denotes the action of $A \in G$ on
$x \in X.$ Suppose for the beginning that the action of $G$ on $X$
is effective, i.e. $A$ acts identically on $X$ only if $A=I$. Then
we claim that both $A(x,\lambda, \zeta)$ and $B^{\rm
T}(x,\lambda,\zeta)$ are Lax matrices for $R$. The claim about $B$
is equivalent to saying that $B(x,\lambda,\zeta)$ is a Lax matrix
for $R_{21}.$

The following argument is illustrated by either the standard or
the ``cubic'' diagram for the Yang-Baxter relation
(Figs. \ref{Fig:YB0},\ref{Fig:YB}). Look at the
values of $z_{12}$ produced by the both parts of the Yang-Baxter
relation (\ref{sYB}): the left-hand side gives
$z_{12}=A(y_1,\mu,\nu)A(x_2,\lambda,\nu)[z]$, while the right-hand
side gives $z_{12}=A(x,\lambda,\nu) A(y,\mu,\nu)[z]$. Now since we
assume that the action of $G$ is effective, we immediately arrive
at the relation
\[
A(x,\lambda,\nu) A(y,\mu,\nu)=A(y_1,\mu,\nu)A(x_2,\lambda,\nu),
\]
which holds whenever $(x_2,y_1)=R(\lambda,\mu)(x,y)$. This
coincides with (\ref{Lax}), an arbitrary parameter $\nu$ playing
the role of the spectral parameter $\zeta$.

Similarly, one could look at the values of $x_{23}$ produced by
the both parts of (\ref{sYB}): the left-hand side gives
$x_{23}=B(z,\nu,\lambda)B(y,\mu,\lambda)[x]$, while the right-hand
side gives $x_{23}=B(y_3,\mu,\lambda) B(z_2,\nu,\lambda)[x]$.
Effectiveness of the action of $G$ again implies:
\[
B(z,\nu,\lambda)B(y,\mu,\lambda)=B(y_3,\mu,\lambda)
B(z_2,\nu,\lambda),
\]
whenever $(y_3,z_2)=R(\mu,\nu)(y,z)$. This turns into (\ref{Lax})
for the transposed matrices $B^{\rm T}$ (or for the inverse
matrices $B^{-1}$); the role of spectral parameter is here played
by an arbitrary parameter $\lambda$.

It should be mentioned that this kind of arguments was first used
to derive Lax representations for 3D consistent discrete equations
on quad-graphs with fields on vertices in \cite{BS, N}. In fact
the 3D consistency condition is the exact analog of the
Yang-Baxter relation for the problems with fields on vertices (see
\cite{ABS}).

In order to cover all the known examples we have to extend the
proposed scheme in the following way. Let us say that $A(x,
\lambda,\zeta)$ gives a {\it projective Lax representation} for
the Yang-Baxter map $R$ if the relation (\ref{Lax}) holds up to
multiplication by a scalar matrix $c I$, where $c$ may depend on
all the variables in the relation. One can easily modify the
arguments from \cite{V} to produce the integrals for the
transfer-maps using the projective Lax matrix: all the ratios of
the eigenvalues of the monodromy matrix are obviously preserved by
these maps.

Assume now that the action of $G = GL_N$ on $X$ is projective,
i.e. scalar matrices are acting trivially and moreover if the
action of $A$ on $X$ is trivial then $A$ is a scalar. Then our
previous considerations show that the matrices
$A(x,\lambda,\zeta)$ and $B^{\rm T}(x,\lambda,\zeta)$ give
projective Lax representations for the corresponding Yang-Baxter
maps (\ref{map}). In practice for a natural choice of matrices
$A$, $B$ in (\ref{map}) we have actually proper Lax
representations, as the following examples show.

\subsection*{Example 1}
{\it Adler's map}.
\medskip

Here $X = {\bf CP}^1$ and the map has the form
\begin{equation}\label{Adler}
\tilde{x} = y-\frac{\lambda-\mu}{x+y}, \qquad
 \tilde{y} = x-\frac{\mu-\lambda}{x+y}.
\end{equation}
This map (modulo additional permutation) first appeared in Adler's
paper \cite{A} as a symmetry of the periodic dressing chain \cite
{VS}. The Lax pair for this map was known from the very beginning
since it comes from re-factorization problem for the matrix
 $$ A(x,\lambda,\zeta) = \Bigl(\begin{array}{cc}
 x &  x^2 + \lambda - \zeta\\ 1 & x \end{array} \Bigl). $$
Our point is that we can actually see this matrix directly in the
map: $$\tilde{y} = x - \frac{\mu - \lambda}{x+y} = \frac{x^2 + xy
-(\mu - \lambda)}{x+y} = A(x, \lambda, \mu)[y],$$ where the group
$G = GL_2$ is acting on ${\bf CP}^1$ by M\"obius transformations.
In this example $B(x,\lambda,\zeta)=A(x,\lambda,\zeta)$, which
reflects the symmetry of the map: $R_{21}=R$.

\subsection*{Example 2}
{\it Interaction of matrix solitons.}
\medskip

One-soliton solutions of the matrix KdV equation $$
U_t+3UU_x+3U_xU+U_{xxx}=0 $$ have the form \cite{G} $$U = 2
\lambda^2 P \sech^2 (\lambda x- 4\lambda^3 t),$$ where the matrix
amplitude $P$ must be a projector: $P^2 = P$, and $\lambda$ is the
parameter measuring the soliton velocity. If we assume that $P$
has rank 1 then $P=\dfrac{\xi\otimes \eta}
{\langle\xi,\eta\rangle}$. Here $\xi$ is a vector in a vector
space $V$ of dimension $N$, $\eta$ is a (co)vector from the dual
space $V^*$, and bracket $\langle\xi,\eta\rangle$ means the
canonical pairing between $V$ and $V^*.$

The change of the matrix amplitudes $P$ of two solitons with the
velocities $\lambda_1$ and $\lambda_2$ after their interaction is
described by the following Yang-Baxter map \cite{G, GV}:
$$R(\lambda_1, \lambda_2): (\xi_1, \eta_1; \xi_2, \eta_2)
\rightarrow (\tilde{\xi}_1, \tilde{\eta}_1; \tilde{\xi}_2,
\tilde{\eta}_2),$$
\begin{equation} \label{fi}
 \tilde{\xi}_1 = \xi_1+\dfrac{2\lambda_2\langle\xi_1,\eta_2\rangle}
 {(\lambda_1-\lambda_2)\langle\xi_2,\eta_2\rangle}\xi_2, \qquad
 \tilde{\eta}_1 = \eta_{1}+\dfrac{2\lambda_2\langle\xi_2,\eta_1\rangle}
 {(\lambda_1-\lambda_2)\langle\xi_2,\eta_2\rangle}\eta_2,
\end{equation}
\begin{equation} \label{fj}
 \tilde{\xi}_2 = \xi_2+\dfrac{2\lambda_1\langle\xi_2,\eta_1\rangle}
 {(\lambda_2-\lambda_1)\langle\xi_1,\eta_1\rangle}\xi_1,\qquad
 \tilde{\eta}_2 = \eta_2+\dfrac{2\lambda_1\langle\xi_1,\eta_2\rangle}
 {(\lambda_2-\lambda_1)\langle\xi_1,\eta_1\rangle}\eta_1.
\end{equation}

In this example $X$ is the set of projectors $P$ of rank 1 which
is the variety ${\bf CP}^{N-1} \times {\bf CP}^{N-1}$, and the
group $G = GL_N$ is acting on the projectors by conjugation (which
corresponds to the natural action of $GL(V)$ on $V \otimes V^*$).

It is easy to see that the formulas (\ref{fi}), (\ref{fj}) are of
the form (\ref{map}) with the matrices $$A(P,\lambda,\zeta)
=B(P,\lambda,\zeta)=I+\frac{2\lambda}{\zeta-\lambda}P=I+ \frac{2
\lambda}{\zeta -\lambda}\cdot\dfrac{\xi\otimes \eta}
{\langle\xi,\eta\rangle}$$ (note that again $R_{21}=R$). Our
results show that the matrix $A(P,\lambda,\zeta)$ gives a
projective Lax representation for the interaction map. In
\cite{GV} it is shown that this is actually a genuine Lax
representation. One can explain in the same way the Lax matrices
for more general Yang-Baxter maps on Grassmannians from \cite{GV}.

\subsection*{Example 3}
{\it Yang-Baxter maps arising from geometric crystals} \cite{Et},
\cite{NY}.
\medskip

Let $X={\bf C}^n$, and define $R:X\times X\rightarrow X\times X$
by the formulas
\begin{equation}\label{E}
 \tilde{x}_j=x_j\,\frac{P_j}{P_{j-1}},\qquad
 \tilde{y}_j=y_j\,\frac{P_{j-1}}{P_j},\qquad j=1,\ldots,n,
\end{equation}
 where
\begin{equation}\label{EP}
 P_j=\sum_{a=1}^n\left(\prod_{k=1}^{a-1}x_{j+k}\prod_{k=a+1}^ny_{j+k}\right)
\end{equation}
(in this formula subscripts $j+k$ are taken (mod $n$)).
 Clearly, the map (\ref{E}) keeps the following subsets invariant:
 $X_\lambda\times X_\mu\subset X\times X$, where
 $X_\lambda=\{(x_1,\ldots,x_n)\in X: \prod_{k=1}^nx_k=\lambda\}$.
 It can be shown that the restriction of $R$ to $X_\lambda\times X_\mu$
 may be written in the form (\ref{map}). For this, the following trick is
 used. Embed this set into
 ${\bf CP}^{n-1}\times{\bf CP}^{n-1}$:
\[
J(x,y)=(z(x),w(y)),\quad z(x)=(1:z_1:\ldots:z_{n-1}), \quad
w(y)=(w_1:\ldots:w_{n-1}:1),
\]
\[
z_j=\prod_{k=1}^{j}x_k\,,\qquad w_j=\prod_{k=j+1}^ny_k\,.
\]
Then it is easy to see that in coordinates $(z,w)$ the map $R$ is
written as
\[
\tilde{z}=B(y,\mu,\lambda)[z]\,,\qquad
\tilde{w}=A(x,\lambda,\mu)[w]\,,
\]
with certain matrices $B,A$ from $G=GL_n$, where the standard
projective action of $GL_n$ on ${\bf CP}^{n-1}$ is used. Moreover,
a simple calculation shows that the inverse matrices are cyclic
two-diagonal:
\begin{equation}\label{B}
B^{-1}(y,\mu,\lambda)=\left(\begin{array}{cccccc}
 y_1 & -1 & 0 & \ldots & 0 & 0\\
 0 & y_2 & -1 & \ldots & 0 & 0\\
 0 & 0 & y_3 & \ldots & 0 & 0 \\
 & \ldots & &&\ldots &\\
 0 & 0 & 0 & \ldots & y_{n-1} & -1\\
 -\lambda &  0 & 0 & \ldots & 0 & y_n
 \end{array}\right),
\end{equation}
\begin{equation}\label{A}
A^{-1}(x,\lambda,\mu)=\left(\begin{array}{cccccc}
 x_1 & 0 & 0 & \ldots & 0 & -\mu\\
 -1 & x_2 & 0 & \ldots & 0 & 0\\
 0 & -1 & x_3 & \ldots & 0 & 0 \\
 & \ldots & &&\ldots &\\
 0 & 0 & 0 & \ldots & x_{n-1} & 0\\
 0 &  0 & 0 & \ldots & -1 & x_n
 \end{array}\right).
\end{equation}
To be more precise the matrices $A,B$ are defined only up to
multiplication by scalar matrices. These scalar matrices are
chosen in (\ref{B}), (\ref{A}) in such a way that the dependence
of the matrices $B^{-1}$, $A^{-1}$ on their ``own'' parameters
($\mu$ and $\lambda$, resp.) drops out, so that the only parameter
remaining in the Lax representation is the spectral one. In other
words, the Lax representation does not depend on the subset
$X_\lambda\times X_\mu$ to which we restricted the map. Note also
that we get this time only {\it one} Lax representation for $R$,
since the matrices $B^{\rm T}$ coincide with $A$. It can be
checked that this is actually a genuine (not only projective) Lax
representation.

As the last remark we would like to mention that our Lax
representation is closely related to the notion of the {\it
structure group} $G_R$ of the Yang-Baxter map $R$ \cite{Et}. It
was shown by Etingof in \cite{Et} that for the map (\ref{E}) the
so-called reduced structure groups $G_R^+$ and $G_R^-$ can be
realized as the subgroups of the loop group $PGL_n(C(\lambda))$
generated by the matrix functions $A^{-1}(x,\cdot,\lambda)$ with $x\in X$,
resp. by $B^{-1}(x,\cdot,\lambda)$ with $x\in X$.

\subsection*{Acknowledgements}

We are grateful to the organizers of the SIDE-V conference in
Giens (21-26 June 2002) where this work was done.

\end{document}